\documentclass[7pt]{article}
\usepackage{amsfonts,amssymb,amsmath,comment}
\usepackage{graphicx,graphics,color}
 \makeatletter
 \@addtoreset{equation}{section}
 \makeatother

 \newcounter{enunciato}[section]

 \newtheorem{ittheorem}{Theorem}
 \newtheorem{itlemma}{Lemma}
 \newtheorem{itproposition}{Proposition}
 \newtheorem{itdefinition}{Definition}
 \newtheorem{itcorollary}{Corollary}
 \newtheorem{itconjecture}{Conjecture}

 \newenvironment{theorem}{\addtocounter{enunciato}{1}
 \begin{ittheorem}}{\end{ittheorem}}

 \newenvironment{lemma}{\addtocounter{enunciato}{1}
 \begin{itlemma}}{\end{itlemma}}

 \newenvironment{corollary}{\addtocounter{enunciato}{1}
 \begin{itcorollary}}{\end{itcorollary}}

\newcommand{\halmos}{\rule{1ex}{1.4ex}}

\def \ba {\begin{array}}
\def \ea {\end{array}}

\def \R {{\mathbb R}}

\def \N {{\mathbb N}}

\def \P {{\mathbb P}}

\def \Leb {{\mathcal{L}}}

\def \ra {\rightarrow}

\def \cp {\mathrm{cap}\,}
\def\O{\Omega}

\date{22 October 2016}
\begin{document}
\title{On P\'olya's inequality for torsional rigidity and first Dirichlet eigenvalue}

\author{M. van den Berg \thanks{MvdB acknowledges support by The Leverhulme Trust
through International Network Grant \emph{Laplacians, Random Walks, Bose Gas,
Quantum Spin Systems.}}, V. Ferone, C. Nitsch, C. Trombetti
 \\
\\
School of Mathematics, University of Bristol\\
University Walk, Bristol BS8 1TW, UK\\
\texttt{mamvdb@bristol.ac.uk}\\ \\
Dipartimento di Matematica e Applicazioni ``Renato Cacciopoli''\\
 Universit\`a degli Studi di Napoli Federico II\\
Via Cintia, Monte S. Angelo, I-80126 Napoli, Italy\\
\texttt{vincenzo.ferone@unina.it}\\
\texttt{c.nitsch@unina.it}\\
\texttt{cristina.trombetti@unina.it}}

\maketitle

\begin{abstract}
Let $\Omega$ be an open set in Euclidean space with finite Lebesgue measure $\vert\Omega\vert$. We
obtain some properties of the set function $F:\Omega\mapsto \R^+$ defined by
$$ F(\Omega)=\frac{T(\Omega)\lambda_1(\Omega)}{\vert\Omega\vert} ,$$
where $T(\Omega)$ and $\lambda_1(\Omega)$ are the torsional rigidity and the first eigenvalue of the Dirichlet Laplacian respectively. We improve the classical P\'olya bound $F(\Omega)\le 1,$ and show that $$F(\Omega)\le 1- \nu_m T(\Omega)|\Omega|^{-1-\frac2m},$$ where $\nu_m$ depends only on $m$. For any $m=2,3,\dots$ and $\epsilon\in (0,1)$ we construct an open set $\Omega_{\epsilon}\subset \R^m$ such that $F(\Omega_{\epsilon})\ge 1-\epsilon$.

\vskip 0.5truecm
\noindent
{\it AMS} 2000 {\it subject classifications.} 49J45, 49R05, 35P15, 47A75, 35J25.\\
{\it Key words and phrases.} Torsional rigidity, first Dirichlet eigenvalue

%\medskip\noindent

\end{abstract}

%%%%%%%%%%%%%%%%%%%%%%%%%%%%%%%%%%%%%%%

\section{Introduction}
\label{Introduction}
Let $\Omega$ be an open set in $\R^m$ with finite Lebesgue measure $\vert\Omega\vert$, and let $v_{\Omega}:\Omega\mapsto \R^+$ denote the corresponding torsion function, i.e. the unique solution of
\begin{equation}\label{e1}
-\Delta v=1,\, v\in H_0^1(\Omega).
\end{equation}
The torsional rigidity of $\Omega$ is defined by $T(\Omega)=\int_{\O}v_\Omega$. As $v_{\Omega}\ge 0$, the torsional rigidity is the $\Leb^1(\O)$ norm of $v_{\Omega}$. The following variational characterisation is well known
\begin{equation*}%\label{e2}
T(\Omega)=  \sup_{w\in H_0^1(\Omega)\setminus\{0\}}\frac{\left(\displaystyle\int_\Omega w\right)^2}{\displaystyle\int_\Omega |Dw|^2}.
\end{equation*}
The torsional rigidity plays a key role in different parts of analysis.
For example the torsional rigidity of a cross section of a beam appears in the computation
of the angular change when a beam of a given length and a given modulus of
rigidity is exposed to a twisting moment \cite{Bandle}, \cite{PSZ}. It also arises
in the calculation of the heat content of sets with time-dependent boundary
conditions \cite{MvdB}, in the definition of gamma convergence \cite{BB}, and
in the study of minimal submanifolds \cite{MP}. Moreover, $T(\O)/|\O|$ equals the expected lifetime of Brownian motion in $\O$ when averaged with respect to the uniform distribution over all starting points $x\in\O$.

Since $\O$ has finite Lebesgue measure the Dirichlet Laplacian acting in $\Leb^2(\O)$ has compact resolvent. We denote the eigenvalues and a corresponding orthonormal basis by $\lambda_1(\O)\le \lambda_2(\O)\le \dots$, and $\{\varphi_1,\varphi_2,\dots\}$ respectively. Recall the following variational characterisation.
\begin{equation*}%\label{e3}
\lambda_1(\Omega)=\inf_{z\in H_0^1(\Omega)\setminus\{0\}}\frac{\displaystyle\int_\Omega|Dz|^2}{\displaystyle\int_\Omega z^2}.
\end{equation*}
A classical inequality of P\'olya, \cite{PSZ}, asserts that the function $F$ defined by
\begin{equation}\label{e4}
F(\Omega)=\frac{T(\Omega)\lambda_1(\Omega)}{\vert\Omega\vert}
\end{equation}
satisfies
\begin{equation}\label{e5}
F(\Omega)\le 1.
\end{equation}
We note that $F$ is scale independent i.e. for any homothety $\alpha\O,\alpha>0$, of $\O$ we have that $F(\alpha\O)=F(\O)$.

The main results of this paper are the following.
\begin{theorem}\label{the4}
 For any open set $\Omega$ with finite Lebesgue measure
\begin{equation*}%\label{a1}
F(\Omega)\le 1- \frac{2m \omega_m^{2/m}}{m+2}\frac{T(\Omega)}{|\Omega|^{1+\frac2m}},
\end{equation*}
where $\omega_m$ is the measure of the ball with radius $1$ in $\R^m$.
\end{theorem}

\begin{theorem}\label{the1}
Let $m=2,3,\dots$. For every $\epsilon>0$ there exists an open connected set $\Omega_{\epsilon}\subset \R^m$ depending on $\epsilon$ such that
\begin{equation}\label{e6}
F(\Omega_{\epsilon})\ge 1-\epsilon.
\end{equation}
\end{theorem}

\begin{corollary}
The variational problem $$\sup\{F(\Omega): \Omega\ \textup{open in}\,  \R^m,\, \vert\Omega\vert=1\}$$ does not have a maximiser.
\end{corollary}
The proof of this corollary is immediate. Indeed, by Theorem \ref{the1}, and P\'olya's inequality the above supremum equals $1$.  Suppose there exists an open set $\Omega$ with $F(\Omega)=1$ and $|\Omega|=1$.
Then $\Omega$ has strictly positive torsional rigidity, and $F(\Omega)<1$ by Theorem \ref{the4} which is a contradiction.

It was shown in \cite[Remark 2.4]{MvdBBV} that
\begin{equation}\label{e7}
\inf\{F(\Omega): \Omega\, \textup{open in}\, \R^m,\, \vert\Omega\vert=1\}=0.
\end{equation}
However, for the restriction of $F(\cdot)$ to the class of convex sets in $\R^m$, we have the following.
\begin{theorem}\label{the2}
\begin{itemize}
\item[\textup{(i)}]
\begin{equation}\label{e8}
\inf\{F(\Omega): \Omega\, \textup{open, convex in\,} \R^m,\, \vert\Omega\vert=1\}\ge \frac{\pi^2}{4m^{m+2}(m+2)}.
\end{equation}
\item[\textup{(ii)}]
\begin{equation}\label{e9}
\inf\{F(\Omega): \Omega\, \, \textup{open, convex in\,} \R^2,\, \vert\Omega\vert=1\}\geq\frac{\pi^2}{48}.
\end{equation}
\end{itemize}
\end{theorem}

Theorem \ref{the1} disproves the conjecture in \cite{MvdBBV} that $F(\Omega)\le \frac{\pi^2}{12}$. Theorem \ref{the3} below goes some way towards proving the $\pi^2/12$ bound for open bounded, planar, convex sets. In order to state our main result for convex sets, we introduce the following notation. For a convex set with finite measure, we denote by $w$ the minimum width of $\Omega$ (or simply the width of $\Omega$), which is obtained by minimising among all pairs of parallel supporting hyperplanes of $\Omega$ the distance between such hyperplanes. The projection of $\Omega$ onto one of the minimising hyperplanes is denoted by $E$. The first eigenvalue of the $(m-1)$-dimensional Dirichlet Laplacian acting in $\Leb^2(E)$ is denoted by $\Lambda$.
\begin{theorem}\label{the3}
\begin{itemize}
\item[\textup{(i)}] If $\Omega$ is an open, bounded, convex set in $\R^m$ with $w$ and $\Lambda$ as above, then
\begin{equation}\label{a2}
F(\Omega)\le \frac{\pi^2}{12}\bigg(1+\frac{3c}{2}+\frac{3c^2}{4}+\frac{c^3}{8}\bigg),
\end{equation}
where
\begin{equation}\label{c22}
c=\bigg(\frac{32w^2\Lambda}{\pi^2}\bigg)^{1/3}.
\end{equation}
\item[\textup{(ii)}]  If $\Omega$ is an open, bounded, convex set in $\R^2$, then
\begin{equation}\label{a3}
F(\Omega)\le 1-\frac{1}{11560}.
\end{equation}
\end{itemize}
\end{theorem}

\begin{corollary}\label{cor1} If $(\O_n)$ is a sequence of bounded convex sets with corresponding sequences $(w_n)$ and $(\Lambda_n)$ such that $\lim_{n\rightarrow \infty}w_n^2\Lambda_n=0$, then $\limsup_{n\rightarrow \infty}F(\O_n)\le \frac{\pi^2}{12}.$
 \end{corollary}

The main idea in the proof of Theorem \ref{the1} is that if $\Omega$ is an open, bounded and connected set, then we can find $x_0\in \Omega$ and $\delta>0$ such that punching a hole in $\Omega$ centered at $x_0$ with radius $\delta$ increases $F$. In the proof of Theorem \ref{the1}, we take an $m$-dimensional cube with side-length $L$ and punch $N^m$ holes with the same radius $\delta$ in a periodic arrangement. We show that we can find $L, N, \delta$ depending on $\epsilon$ (and $m$) such that the corresponding value of $F$ for the punched cube exceeds $1-\epsilon$. As mentioned above, $F$ is invariant under homotheties, and so we could have chosen $L=1$. However, it is convenient to keep $L$ undetermined so that we have a homothety or scaling check in the various bounds.

To see that punching a hole increases $F$, we take $\Omega$ open, bounded, connected, and with smooth boundary. Let $\varphi_1\in H_0^1(\Omega)$ be a Dirichlet eigenfunction corresponding to $\lambda_1(\Omega)$, and let $v_{\Omega}$ be the solution of \eqref{e1}. We observe that $$\lambda_1(\Omega)< \frac{\|Dv_{\Omega}\|^2_{\Leb^2(\O)}}{\|v_{\Omega}\|^2_{\Leb^2(\O)}},$$
implies
\begin{equation*}%\label{e14}
\int_\Omega \left(T(\Omega)\varphi_1^2-\lambda_1(\Omega)v_{\Omega}^2\right)=\|v_{\Omega}\|^2_{\Leb^2(\O)}\left(\frac{1}{\|v_{\Omega}\|^2_{\Leb^2(\O)}}\displaystyle\int_\Omega v_{\Omega}-\lambda_1(\Omega)\right)=\|v_{\Omega}\|^2_{\Leb^2(\O)}\left(\frac{\|Dv_{\Omega}\|^2_{\Leb^2(\O)}}{\|v_{\Omega}\|^2_{\Leb^2(\O)}}-\lambda_1(\Omega)\right)>0.
\end{equation*}
So there exists $x_0\in\Omega$ such that
\begin{equation*}%\label{e13}
T(\Omega)\varphi_1^2(x_0)-\lambda_1(\Omega)v_{\Omega}^2(x_0)>0.
\end{equation*}
Let $\Omega_{\delta,x_0}= \Omega\setminus B(x_0;\delta),$ where $B(x_0;\delta)$ is the closed ball of radius $\delta>0$ centered at $x_0$.
We want to show that if $\delta$ is small enough, then $F(\Omega_{\delta,x_0})>F(\Omega)$. In the planar case $m=2$, a classical asymptotic formula (see, for instance, \cite[Theorem 1.4.1]{HEN} and the references therein) gives that
\begin{equation}\label{e10}\lambda_1(\Omega_{\delta,x_0})=\lambda_1(\Omega)+\frac{2\pi}{-\log\delta}
\varphi_1^2(x_0)+o\left(\frac1{|\log\delta|}\right),\,\delta\downarrow 0.
 \end{equation}
Moreover, from \cite[Theorem 8.1.6]{MNP}, we have that
\begin{equation}\label{e11}T(\Omega_{\delta,x_0})=T(\Omega)-\frac{2\pi}{-\log\delta}{v_{\Omega}^2(x_0)}+o\left(\frac1{|\log\delta|}\right),\, \delta\downarrow 0.
\end{equation}
By \eqref{e10} and \eqref{e11}, we have that
\begin{equation*}%\label{e12}
\frac{T(\Omega_{\delta,x_0})\lambda_1(\Omega_{\delta,x_0})}{|\Omega_{\delta,x_0}|}
=\frac{T(\Omega)\lambda_1(\Omega)}{|\Omega|}+\frac{2\pi}{\left(-\log\delta\right)|\Omega|}\left(T(\Omega)\varphi_1^2(x_0)
-\lambda_1(\Omega)v_{\Omega}^2(x_0)\right)+o\left(\frac1{|\log\delta|}\right), \, \delta\downarrow 0.
\end{equation*}
Hence $F(\Omega_{\delta,x_0})>F(\Omega)$ for $\delta$ sufficiently small. The same calculation works in the higher dimensional setting ($m>2$) replacing $\frac{2\pi}{-\log \delta}$ by the Newtonian capacity of $B(x_0;\delta)$ in \eqref{e10} and \eqref{e11} (see for example \cite[Theorem 1.4.1]{HEN} and \cite[Theorem 8.1.4]{MNP}, respectively).

%%%%%%%%%%%%%%%%%%%%%%%%%%%%%%%%%%%%%%%%%%%%%%%%%%%%%%%%%%%%%%%%%%%%%%%%%%%%
This paper is organised as follows. In Section \ref{sec4} we prove Theorem \ref{the4}. In Section \ref{sec2} we prove Theorem \ref{the1}, and in Section \ref{sec3} we prove Theorems \ref{the2} and \ref{the3} respectively.
\section{Proof of Theorem \ref{the4}}\label{sec4}

Let $v_{\Omega}$ be the torsion function of $\Omega$. By choosing $v_{\Omega}$ as a test function for the Rayleigh quotient for $\lambda_1(\Omega)$, we obtain that
\begin{equation*}%\label{a4}
\lambda_1(\Omega) \le \frac{\displaystyle\int_{\Omega} v_{\Omega}}{\displaystyle\int_{\Omega} v_{\Omega}^2}.
\end{equation*}
Hence
\begin{equation*}%\label{a5}
F(\Omega) \le \frac {\left(\displaystyle\int_\Omega v_{\Omega}\right)^2}{\left(\displaystyle\int_\Omega v_{\Omega}^2\right) |\Omega|}.
\end{equation*}
Let $M= \sup_{\Omega} v_{\Omega}$. For $\theta \in [0,M]$, we define
\begin{equation*}%\label{a6}
 \mu (\theta) = |\{ x \in \Omega :v_{\Omega}(x)>\theta\}|.
\end{equation*}
We have that
\begin{equation*}%\label{a7}
\int_{\Omega}v_{\Omega}=\int_0^M \mu(\theta)\, d\theta,
\end{equation*}
and
\begin{equation*}%\label{a8}
\int_{\Omega}v_{\Omega}^2=\int_0^M 2\theta\mu(\theta)\, d\theta.
\end{equation*}
For every $\theta \in (0,M)$, we have that
\begin{equation}\label{a9}\mu(\theta) \le (|\Omega|^{2/m} - 2m\omega_m^{2/m} \theta)^{m/2}.
\end{equation}
Indeed, since $v_{\Omega}$ satisfies the torsion equation \eqref{e1} in $\Omega$, arguing similarly to \cite{Talenti}, we have that for $\theta \in (0,M)$,
\begin{equation}\label{a10}
\mu(\theta) = \int_{\{v_{\Omega}=\theta\}} |Dv_{\Omega}| \, d \mathcal{H}^{m-1},
\end{equation}
and
\begin{equation*}%\label{a11}
-\mu'(\theta) \ge \int_{\{v_{\Omega}=\theta\}} \frac{1}{|Dv_{\Omega}|}\, d \mathcal{H}^{m-1}.
\end{equation*}
Denote the perimeter of a measurable set $A$ by $\textup{Per}(A)$.
Applying H\"older's  inequality to $\textup{Per}(\{v_{\Omega}>\theta\})= \int_{\{v=\theta\}} \, d \mathcal{H}^{m-1}$,  we obtain that
\begin{equation}\label{a12}
\textup{Per}(\{v_{\Omega}>\theta\})^2  \leq \mu(\theta) (-\mu'(\theta)).
\end{equation}
By the isoperimetric inequality we have that
\begin{equation*}%\label{a13}
\textup{Per}(\{v_{\Omega}>\theta\})\ge m\omega_m^{1/m} \mu(\theta)^{(m-1)/m}.
\end{equation*}
This, together with \eqref{a12}, gives the differential inequality
\begin{equation*}%\label{a14}
 m^2\omega_m^{2/m} \le - \mu(\theta)^{\frac{2}{m}-1}\mu'(\theta).
\end{equation*}
Integrating this differential inequality gives \eqref{a9}.

For $t\in [0,M]$, define
\begin{equation}\label{a15}
Q(t) = \left( \displaystyle\int_0^t \mu(\theta) \,d\theta \right)^2  - 2\left(\displaystyle\int_{0}^t \theta \mu(\theta) \, d\theta \right) |\Omega|.
\end{equation}
Using \eqref{a9} and \eqref{a15}, it is straightforward to verify that
\begin{equation*}%\label{a16}
 Q'(t) \leq  \frac{|\Omega|^{\frac{m+2}{m}}}{m(m+2)\omega_m^{2/m}}
 \left[ 1- \left(1- \frac{2m\omega_m^{2/m}}{|\Omega|^{2/m}} t\right)^{\frac{m+2}{2}}\right] 2 \mu(t) -2t \mu(t) |\Omega|.
 \end{equation*}
The inequality $(1+y)^{\alpha}\ge 1+\alpha y+ {y^2}, \alpha\ge 2, \, y \ge -1$ then gives that
 \begin{equation}\label{a17}
 Q'(t) \leq - |\Omega|^{1-\frac{2}{m}}  \frac{8m \omega_m^{2/m}}{m+2} {\mu(t)} t^2.
 \end{equation}
Integrating \eqref{a17} over $[0,M]$ and using the fact that $Q(0)=0$ gives that
\begin{equation*}%\label{a18}
Q(M)\le -|\Omega|^{1-\frac{2}{m}}  \frac{8m \omega_m^{2/m}}{m+2} \int_0^M  {\mu(t)} t^2\, dt.
\end{equation*}
H\"older's inequality then yields that
\begin{equation*}%\label{a19}
Q(M) \leq -|\Omega|^{1-\frac{2}{m}}  \frac{2m \omega_m^{2/m}}{m+2} \frac{ \left(\int_0^M 2t \mu(t) dt \right)^2} {\int_0^M \mu(t) dt} = -
 \frac{2m \omega_m^{2/m}}{m+2}   |\Omega|^{1-\frac{2}{m}}  \frac{ \left(\displaystyle\int_\Omega v_\Omega^2\right)^2} {\displaystyle\int_\Omega v_\Omega}.
\end{equation*}
Using the expression for $Q$ and H\"older's inequality gives that
\begin{equation*}%\label{a20}
 \left[\frac {\left(\displaystyle\int_\Omega v_\Omega\right)^2}{\left(\displaystyle\int_\Omega v_\Omega^2\right) |\Omega|} -1\right]  \leq
-  \frac{2m \omega_m^{2/m}}{m+2} \frac{T(\Omega)}{ |\Omega|^{1+\frac{2}{m}}}.
\end{equation*}
This concludes the proof of Theorem \ref{the4}.
 \hspace*{\fill }$\square $

\section{Proof of Theorem \ref{the1}}\label{sec2}
 In this section we provide an example of an open connected set $\Omega_{\epsilon}$ in $\R^m$ which satisfies \eqref{e6}. As the technical tools depend heavily on the relation between torsional rigidity and heat equation we recall some of the essential ingredients in Section \ref{3.1} below. The necessary bounds for the first eigenfunction and eigenvalue with Dirichlet boundary conditions on a ball centred in an $m$-dimensional cube with Neumann boundary conditions will be obtained in Section \ref{3.2}. The proof of Theorem \ref{the1} will be deferred to Section \ref{3.3}.

\subsection{Heat equation and torsional rigidity}\label{3.1}
We denote the Dirichlet
heat kernel for $\O$ by $p_{\Omega}(x,y;t)$, $x,y\in\Omega$, $t>0$. The integral defined by
\begin{equation*}%\label{e15}
u_\O(x;t) = \int_\O dy\, p_{\Omega}(x,y;t)
\end{equation*}
is the solution of
\begin{equation}\label{e16a}
\frac{\partial u(x;t)}{\partial t} = \Delta u(x;t),\, x\in\O,\, t>0,
\end{equation}
\begin{equation}\label{e16b}
\lim_{t\downarrow 0}u(\cdot;t) = 1\, \textup{in}\,\Leb^1(\O),
\end{equation}
\begin{equation}\label{e16c}
u(\cdot;t)\in  H_0^1(\O),\, t>0.
\end{equation}
The interpretation of \eqref{e16a}, \eqref{e16b}, and \eqref{e16c} is that $u_\O(x;t)$ represents the temperature at point $x$ at time $t$
when the initial temperature in $\O$ is $1$ and the temperature of $\partial\O$ is
$0$ for all $t>0$. The heat content of $\Omega$ at time $t$ is defined as
\begin{equation*}%\label{e17}
H_\O(t) = \int_\O u_\O(x;t)\,dx.
\end{equation*}
The Dirichlet heat kernel for $\O$ has the following eigenfunction expansion:
\begin{equation}
\label{e19}
p_\O(x,y;t) = \sum_{j\in\N} e^{-t\lambda_j(\O)} \varphi_j(x)\varphi_j(y).
\end{equation}
It follows from Parseval's formula that
\begin{equation}
\label{e20}
H_\O(t) = \sum_{j\in\N} e^{-t\lambda_j(\O)} \left(\int_\O \varphi_j\right)^2
\leq e^{-t\lambda_1(\O)} \sum_{j\in\N} \left(\int_\O \varphi_j\right)^2
= e^{-t\lambda_1(\O)} |\O|.
\end{equation}
The solution of \eqref{e1} is given by
\begin{equation*}
%\label{e21}
v_\O(x) = \int_0^\infty u_{\Omega}(x;t)\,dt.
\end{equation*}
It follows that
\begin{equation}
\label{e22}
T(\O) = \int_0^\infty H_{\O}(t)\,dt,
\end{equation}
i.e., the torsional rigidity is the integral of the heat content.
 By the first identity in \eqref{e20}, \eqref{e22}, and Fubini's theorem we have that
\begin{align}\label{e23}
T(\O) &= \sum_{j\in\N} \lambda_j(\O)^{-1} \left(\int_\O \varphi_j\right)^2\nonumber \\ &
\le\lambda_1(\O)^{-1}\sum_{j\in\N}  \left(\int_\O \varphi_j\right)^2\nonumber \\ &
=\lambda_1(\O)^{-1}\vert\O\vert,
\end{align}
where we have used Parseval's identity in the last equality above. This implies P\'olya's bound \eqref{e5}. The bound also follows by \eqref{e20} and \eqref{e22}.

 By the first identity in \eqref{e23} we obtain that
\begin{equation*}%\label{e24}
T(\O) \geq  \lambda_1(\O)^{-1} \left(\int_\O \varphi_1\right)^2.
\end{equation*}

 \subsection{Eigenfunction and eigenvalue bounds}\label{3.2}
We introduce the following notation. Let $\Omega_L=(-\frac{L}{2},\frac{L}{2})^m$ be an open cube in $\R^m$ with measure $L^m$, and let $K$ be a compact subset of $\O_L$. We denote the first eigenvalue of the Laplacian acting in $\Leb^2(\O_L- K)$ with Neumann boundary conditions on $\partial \O_L$ and Dirichlet boundary conditions on $\partial K$ by $\mu_{1, K ,L}$. We denote the corresponding normalised eigenfunction by $\varphi_{1, K ,L}$.

The following shows that the $\Leb^1$ norm of the first eigenfunction converges to $L^{m/2}$ as $\mu_{1, K ,L}\downarrow 0$.
\begin{lemma}\label{lem1}
If $m=2,3,4,\dots,$ then
\begin{equation}\label{e252}
L^m\left(1-\left(\frac{4mL^2\mu_{1,K,L}}{3e}\right)^{1/2}\right)\le \|\varphi_{1,K,L}\|_{\Leb^1(\O_L-K)}^2\le L^m.
\end{equation}
\end{lemma}
{\it Proof.}
To prove \eqref{e252}, we note that by Cauchy-Schwarz,
\begin{equation}\label{e38}
\|\varphi_{1,K,L}\|_{\Leb^1(\O_L-K)}^2\le \vert \O_L- K\vert\le \vert \O_L\vert=L^m.
\end{equation}
This proves the right-hand side of \eqref{e252}. To prove the left-hand side of \eqref{e252}, we denote the heat kernel with Neumann boundary conditions on $\partial \O_L$ and Dirichlet boundary conditions
on $\partial K$ by $\pi_{K,L}(x,y;t)$. By the eigenfunction expansion of $\pi_{K,L}(x,y;t)$, we have for $t>0$ that
\begin{align*}%\label{e254}
e^{-t\mu_{1,K,L}}\varphi_{1,K,L}(x)^2&\le \pi_{K,L}(x,x;t)\le \pi_{\O_L}(x,x;t)\nonumber \\ &
\le L^{-m}\left(1+2\sum_{j=1}^{\infty}e^{-t\pi^2j^2/L^2}\right)^{m}\nonumber \\ &\le L^{-m}\left(1+\sum_{j=1}^{\infty}\frac{2L^2}{et\pi^2j^2}\right)^m\nonumber \\ &
=L^{-m}\left(1+\frac{L^2}{3et}\right)^m,
\end{align*}
where $\pi_{\O_L}(x,y;t)$ is the Neumann heat kernel for the cube $\O_L$, and where we have used the eigenfunction expansion of the latter together with separation of variables. Taking the supremum over all $x\in \O_L- K$ gives that
\begin{align*}%\label{e255}
\|\varphi_{1,K,L}\|^2_{\Leb^{\infty}(\O_L-K)}\le e^{t\mu_{1,K,L}}L^{-m}\left(1+\frac{L^2}{3et}\right)^m.
\end{align*}
Furthermore, since  $\|\varphi_{1,K,L}\|^2_{\Leb^2(\O_L-K)}=1$, we have by the positivity of $\varphi_{1,K,L}$ that
\begin{align*}%\label{e256}
\|\varphi_{1,K,L}\|_{\Leb^1(\O_L-K)}^2&\ge \|\varphi_{1,K,L}\|_{\Leb^{\infty}(\O_L-K)}^{-2}\ge L^me^{-t\mu_{1,K,L}}\left(1+\frac{L^2}{3et}\right)^{-m}\nonumber \\ &\ge
L^m\left(1-t\mu_{1,K,L}-\frac{mL^2}{3et}\right).
\end{align*}
We choose $t>0$ as to maximise the right-hand side above. This proves the left-hand side of \eqref{e252}.
\hspace*{\fill }$\square $

In the sequel we need upper and lower bounds for the first Dirichlet eigenvalue $\mu_{1,K,L}$ where $K=B(0;\delta)\subset \O_L.$
These were obtained for general compact sets $K\subset \Omega_L\subset \R^m,\, m=3,4,...$ in \cite{MT1} and \cite{MT2} in terms of the Newtonian capacity $\cp(K)$ of $K$ in $\R^m$. The various $m$-dependent constants in \cite[Propositions 2.2, 2.3, 2.4]{MT1} and in \cite[Theorem A]{MT2} have not been evaluated. We supply these in the Lemmas \ref{lem2} and \ref{lem3} below. We consider general compact subsets as the proofs (for $m=3,...$) are hardly more involved than the special case of a ball.

\begin{lemma}\label{lem2}
\begin{enumerate}
\item[\textup{(i)}]If $m=3,4,\dots$ and if $K\subset \O_L$, then
\begin{equation}\label{e45a}
\mu_{1,K,L}\ge k_m\frac{\cp(K)}{L^{m}},
\end{equation}
where
\begin{equation}\label{e45b}
k_m=\int_0^1ds\, (4\pi s)^{-m/2}e^{-m/(4s)}.
\end{equation}
\item[\textup{(ii)}] If $m=3,4,\dots$ and if $K\subset \O_L$ with $\cp(K)\le \frac{1}{16}L^{m-2},$ then
\begin{equation}\label{e45c}
\mu_{1,K,L}\le 2\pi m\frac{\cp(K)}{L^{m}}.
\end{equation}
\end{enumerate}
\end{lemma}

{\it Proof.}
By the $\Leb^2$-eigenfunction expansion of $\pi_{K,L}(x,y;t)$ we have that
\begin{equation}\label{e45e}
e^{-t\mu_{1,K,L}}\varphi_{1,K,L}(x)=\int_{\O_L-K}dy\,\pi_{K,L}(x,y;t)\varphi_{1,K,L}(y).
\end{equation}
As in \cite{MT1} and \cite{MT2}, we introduce some Brownian motion tools. Let $(\tilde{B}(s),s\ge 0;\tilde{\P}_x,x\in \overline{\O_L})$
be Brownian motion with reflection on $\partial\O_L$. For a compact subset $K\subset \O_L$ we let
\begin{equation}\label{e45f}
\tilde{\tau}_K=\inf\{s\ge 0: \tilde{B}(s)\in K\}.
\end{equation}
Then
\begin{equation}\label{e45g}
\tilde{\P}_x[\tilde{\tau}_{K}>t]=\int_{\O_L-K}dy\, \pi_{K,L}(x,y;t),
\end{equation}
Integrating both sides of \eqref{e45e} with respect to $x$ over $\O_L-K$ gives, with \eqref{e45g}, that
\begin{align*}%\label{e45h}
e^{-t\mu_{1,K,L}}\int_{\O_L-K}dx\, \varphi_{1,K,L}(x)&=\int_{\O_L-K}dy\,\tilde{\P}_y[\tilde{\tau}_K>t]\varphi_{1,K,L}(y)\nonumber \\ &=
\int_{\O_L-K}dx\, \varphi_{1,K,L}(x)-\int_{\O_L-K}dy\,\tilde{\P}_y[\tilde{\tau}_K\le t]\varphi_{1,K,L}(y).
\end{align*}
It follows that
\begin{align*}%\label{e45i}
\mu_{1,K,L}&=-\frac{1}{t}\log\bigg(1-\frac{\int_{\O_L-K}dy\,\tilde{\P}_y[\tilde{\tau}_K\le t]\varphi_{1,K,L}(y)}{\int_{\O_L-K}dy\,\varphi_{1,K,L}(y)}\bigg)\nonumber \\ &
\ge \frac{1}{t}\frac{\int_{\O_L-K}dy\,\tilde{\P}_y[\tilde{\tau}_K\le t]\varphi_{1,K,L}(y)}{\int_{\O_L-K}dy\,\varphi_{1,K,L}(y)}\nonumber \\ &
\ge \frac{1}{t}\inf_{x\in \O_L-K}\tilde{\P}_x[\tilde{\tau}_K\le t].
\end{align*}
Following \cite[p.449]{MT2}, we define $\tilde{K}$ as the subset of $\R^m$ by the method of images, so that in each tiling $L$-cube of $\R^m$ we have a reflected image of $K$.
Let $({B}(s),s\ge 0;{\P}_x,x\in \R^m)$ be Brownian motion on $\R^m$, and define the first hitting time of a closed set $A$ by
\begin{equation}\label{e45j}
\tau_A=\inf\{s\ge 0: B(s)\in A\},
\end{equation}
Then
\begin{equation*}%\label{e45k}
\tilde{\P}_x[\tilde{\tau}_K\le t]=\P_x[\tau_{\tilde{K}}\le t]\ge \P_x[\tau_K\le t].
\end{equation*}
For a compact set $K\subset \R^m$, we define the last exit time by
\begin{equation*}%\label{e45l}
L_K=\sup\{s\ge 0: B(s)\in K\},
\end{equation*}
where we put $L_K=+\infty$ if the supremum is over the empty set.
Then $\P_x[\tau_K\le t]\ge \P_x[L_K\le t]$. By \cite{PS},
we have that
\begin{equation}\label{e45m}
    \mathbb {P}_x[L_K<t]=\int\mu_K(dy)\int^t_0 p(x,y;s)ds,
\end{equation}
where
\begin{equation}\label{e45n}
p(x,y;s)=(4\pi s)^{-m/2}e^{-|x-y|^2/(4s)},
\end{equation}
and where $\mu_K(dy)$ is the equilibrium measure of the compact $K$. Next we choose $t=L^2$. By the above, we have that
\begin{equation}\label{e45o}
\mu_{1,K,L}\ge L^{-2}\inf_{x\in \O_L-K}\int\mu_K(dy)\int^{L^2}_0ds\, p(x,y;s).
\end{equation}
For $y\in K$ and $x\in \O_L-K,$ we have that $|x-y|\le \textup{diam}(\O_L)=mL^2$.
So, by \eqref{e45o}, we conclude that
\begin{equation*}%\label{e45p}
\mu_{1,K,L}\ge L^{-2}\int\mu_K(dy)\int^{L^2}_0ds\, (4\pi s)^{-m/2}e^{-mL^2/(4s)}=k_m\frac{\cp(K)}{L^m},
\end{equation*}
where $k_m$ is given by \eqref{e45b}. This proves part (i) of the lemma.

To prove part (ii) of the lemma, we follow the Remark on p.451 in \cite{MT2}, and define the trial function
\begin{equation}\label{f1}
\psi(x)=1-\kappa_m^{-1}\int \mu_K(dy)\, |x-y|^{2-m},
\end{equation}
where $$\kappa_m=\frac{4\pi^{m/2}}{\Gamma((m-2)/2)},$$
is the Newtonian capacity of the ball with radius $1$ in $\R^m$.
Then
\begin{equation*}%\label{f2}
|D\psi|(x)\le \kappa_m^{-1}(m-2)\int \mu_K(dy)\, |x-y|^{1-m}.
\end{equation*}
Hence
\begin{equation}\label{f3}
\lVert D\psi \rVert _{\Leb^2(\O_L-K)}^2\le \kappa_m^{-2}(m-2)^2\int \mu_K(dy)\,\int \mu_K(dy')\, \int_{\R^m}\,dx\,|x-y|^{1-m}|x-y'|^{1-m}.
\end{equation}
In order to compute the integral with respect to $x$ over $\R^m$, we write
\begin{equation}\label{f4}
|x-y|^{1-m}=\frac{2\pi^{m/2}}{\Gamma((m-1)/2)}\int_0^{\infty}\, \frac{ds}{s^{1/2}}p(x,y;s).
\end{equation}
By Tonelli's theorem, \eqref{f4}, and the semigroup property of the heat kernel, we have that
\begin{align}\label{f5}
\int_{\R^m}\,dx\,|x-y|^{1-m}&|x-y'|^{1-m}=\bigg(\frac{2\pi^{m/2}}{\Gamma((m-1)/2)}\bigg)^2\int_0^{\infty}\,\int_{\R^m}dx\, \frac{ds}{s^{1/2}}p(x,y;s)\int_0^{\infty}\,\frac{ds'}{s'^{1/2}}p(x,y';s')\nonumber \\ & =
\bigg(\frac{2\pi^{m/2}}{\Gamma((m-1)/2)}\bigg)^2\int_0^{\infty}\frac{ds}{s^{1/2}}\int_0^{\infty}\frac{ds'}{s'^{1/2}}p(y,y';s+s').
\end{align}
Changing variables $s=\sigma^2,s'=\sigma'^2$ gives that the right-hand side above equals
\begin{align}\label{f6}
4\bigg(\frac{2\pi^{m/2}}{\Gamma((m-1)/2)}\bigg)^2\int_0^{\infty}\int_0^{\infty}d\sigma\, d\sigma' p(y,y';\sigma^2+\sigma'^2)=\frac{\pi^{(m+2)/2}}{\Gamma((m-1)/2)^2}\Gamma((m-2)/2)|y-y'|^{2-m}.
\end{align}

By \eqref{e45m} and \eqref{e45n}, we have that for $y\in K$,
\begin{align}\label{f7}
1=\P_y[L_K<\infty]=\kappa_m^{-1}\int \mu_K(dy')|y-y'|^{2-m}.
\end{align}
Putting \eqref{f3}-\eqref{f7} together gives that
\begin{align}\label{f8}
\lVert D\psi \rVert _{\Leb^2(\O_L-K)}^2\le \pi\bigg(\frac{\Gamma(m/2)}{\Gamma((m-1)/2)}\bigg)^2\cp(K)\le \pi m \, \cp(K).
\end{align}
The last inequality in \eqref{f8} follows from uniform bounds on the $\Gamma$ function. See for example \cite[6.1.38]{AS}.
We obtain a lower bound for $\lVert \psi \rVert _{\Leb^2(\O_L-K)}$ as follows. By \eqref{f1},
\begin{align*}%\label{f9}
\lVert \psi \rVert _{\Leb^2(\O_L-K)}^2&\ge \int_{\O_L-K}\,dx(1-2\kappa_m^{-1}\int \mu_K(dy)\, |x-y|^{2-m})\nonumber \\ &
=|\O_L-K|- 2\kappa_m^{-1}\int \mu_K(dy)\int_{\O_L-K}dx\,|x-y|^{2-m}.
\end{align*}
By rearrangement, we have that
\begin{equation*}%\label{f10}
\int_{\O_L-K}dx\,|x-y|^{2-m}\le \int_{\O_L^*}dx\,|x|^{2-m}= 2^{-1}m\omega_m^{(m-2)/m}L^{2},
\end{equation*}
where $\O_L^*$ is the ball centered at $0$ with the same measure as $\O_L.$
Hence
\begin{align}\label{f11}
\lVert \psi \rVert _{\Leb^2(\O_L-K)}^2&\ge |\O_L|-|K|-\kappa_m^{-1}m\omega_m^{(m-2)/m}\cp(K)L^{2}\nonumber \\ &
\ge L^m-|K|-\cp(K)L^{2}.
\end{align}
By the classical isoperimetric inequality for the Newtonian capacity of $K$,
\begin{equation}\label{f12}
|K|\le \omega_m\bigg(\frac{\cp(K)}{\kappa_m}\bigg)^{m/(m-2)}\le 7 \cp(K)^{m/(m-2)},
\end{equation}
where, in the last inequality, we have used \cite[6.1.38]{AS} once more.
From \eqref{f11} and \eqref{f12}, we obtain that
\begin{equation}\label{f13}
\lVert \psi \rVert _{\Leb^2(\O_L)}^2\ge L^m-7\cp(K)^{m/(m-2)}-\cp(K)L^{2}.
\end{equation}
If $\cp(K)\le cL^{m-2}$ then the right-hand side of \eqref{f13} is bounded from below by $L^m/2$ provided
$7c^{m/(m-2)}+c\le\frac12$. This clearly holds for all $c\le \frac{1}{16}$. So if $\cp(K)\le \frac{1}{16}L^{m-2}$, then
$\lVert \psi \rVert _{\Leb^2(\O_L-K)}^2\ge L^m/2$. This, together with \eqref{f8}, completes the proof of \eqref{e45c}.
\hspace*{\fill }$\square $

For the two-dimensional case we only consider $K=B(0;\delta)\subset \O_L$.
\begin{lemma}\label{lem3} For $m=2$ and $\delta<\frac{L}{6}$,
\begin{equation}\label{f14}
\frac{1}{100L^2}\bigg(\log \frac{L}{2\delta}\bigg)^{-1}\le \mu_{1,B(0;\delta),L}\le \frac{8\pi}{(4-\pi)L^2}\bigg(\log \frac{L}{2\delta}\bigg)^{-1}.
\end{equation}
\end{lemma}

{\it Proof.}
We define
\begin{equation*}%\label{f15}
\psi(x)= \begin{cases}
\dfrac{\log \frac{|x|}{\delta}}{\log \frac{L}{2\delta}},\, \delta\le |x|\le \frac{L}{2},\\
\\
1,\, x\in \O_L\cap\{|x|>\frac{L}{2}\}.
\end{cases}
\end{equation*}
Then
\begin{equation*}%\label{f16}
\lVert D \psi\rVert_{\Leb^2(\O_L-B(0;\delta))}^2=2\pi \bigg(\log \frac{L}{2\delta}\bigg)^{-1},
\end{equation*}
and
\begin{equation*}%\label{f17}
\lVert \psi\rVert_{\Leb^2(\O_L-B(0;\delta))}^2\ge |\O_L\cap\{|x|>\tfrac{L}{2}\}|=\bigg(1-\frac{\pi}{4}\bigg)L^2.
\end{equation*}
This proves the upper bound in \eqref{f14}.

To prove the lower bound we use the method of descent as in \cite[p.451]{MT2}, and observe that for $m=2$, $\mu_{1,B(0;\delta),L}$ equals the bottom of the spectrum of the
Laplacian with Neumann boundary conditions on the boundary of the cube $\Omega_L=(-\frac{L}{2},\frac{L}{2})^3$, and Dirichlet boundary conditions on the cylinder $C_{L,\delta}=\{(x_1,x_2,x_3):-\frac{L}{2}<x_1<\frac{L}{2}, x_2^2+x_3^2<\delta^2\}$ of height $L$ and radius $\delta$ through the centre of that cube. By the lower bound in Lemma \ref{lem2} for $m=3$, we obtain that
\begin{equation}\label{f18}
\mu_{1,B(0;\delta),L}\ge k_3\frac{\cp(C_{L,\delta})}{L^3}.
\end{equation}
It remains to find a lower bound for $\cp(C_{L,\delta})$. To that end, we follow similar arguments to the proof of \cite[Proposition 3.4, pp.67,68]{PS}.
We consider the $N$ balls $B_1,\dots,B_N$ with radii $\delta$ and centres $(-\frac{L}{2}+(2j-1)\delta,0,0),\, j=1,\dots N$ where $N=\lfloor \frac{L}{2\delta} \rfloor$. Recall that for any compact set $K\subset \R^3$,
\begin{equation*}%\label{f19}
\cp(K)=\sup\bigg\{ \bigg(\iint \frac{\mu_K(dx)\mu_K(dy)}{4\pi|x-y|}\bigg)^{-1}:\, \mu\, \textup{is a probability measure supported on $K$}   \bigg\}.
\end{equation*}
By monotonicity, we have that
$$\cp(C_{L,\delta})\ge \cp(\cup_{j=1}^NB_j). $$
To bound the latter, we let $\sigma_j$ be the surface measure on the boundary of the $j$th ball, and let
\begin{equation*}%\label{f20}
\mu=\frac{1}{4\pi N\delta^2}\sum_{j=1}^N \sigma_j.
\end{equation*}
We wish to find an upper bound for the energy
\begin{equation}\label{f21}
\frac{1}{(4\pi N\delta^2)^2}\sum_{j=1}^N\sum_{k=1}^N\iint\frac{\sigma_j(dx)\sigma_k(dy)}{4\pi |x-y|}.
\end{equation}
If $N=1$ then the expression above equals the inverse of $\cp(B(0;\delta))$. The contribution from the $N$ terms with $j=k$ in \eqref{f21} equals $\frac{1}{4\pi N\delta}$.
Furthermore, the contribution of the terms with $|j-k|=1$ in \eqref{f21} is bounded by  $\frac{N-1}{(4\pi N\delta^2)^2}\iint\frac{\sigma_1(dx)\sigma_2(dy)}{4\pi |x-y|}$.
As $\delta^{-1}d\sigma_j$ is the equilibrium measure for the $j$th ball, we have that $\int\frac{\delta^{-1}\sigma_2(dy)}{4\pi |x-y|}\le 1$. We conclude that
\begin{equation*}%\label{f22}
\frac{N-1}{(4\pi N\delta^2)^2}\iint\frac{\sigma_1(dx)\sigma_2(dy)}{4\pi |x-y|}\le \frac{N-1}{4\pi N^2\delta^3}\int\sigma_1(dx)\le \frac{1}{4\pi N\delta}.
\end{equation*}
Similarly the contribution of the terms with $|j-k|=2$ in \eqref{f21} is bounded by
\begin{equation*}%\label{f23}
\frac{N-2}{(4\pi N\delta^2)^2}\iint\frac{\sigma_1(dx)\sigma_3(dy)}{4\pi |x-y|}\le \frac{1}{4\pi N\delta}.
\end{equation*}
It remains to find an upper bound for the terms in \eqref{f21} for $|j-k|\ge 3$. For $x,y$ on the surface of the $j,k$th balls we have that $|x-y|\ge 2|k-j-1|\delta$.
Hence the contribution from the terms with $|j-k|\ge 3$ in \eqref{f21} is bounded from above by
\begin{align*}%\label{f24}
&\frac{1}{(4\pi N\delta^2)^2}\sum_{k=1}^N\bigg(\sum_{j\ge k+3}^N\frac{1}{8\pi(j-k-1)\delta}+\sum_{j\le k-3}\frac{1}{8\pi(k-1-j)\delta}\bigg)\sigma_1(B_1)\sigma_2(B_2)\nonumber\\ &
=\frac{1}{8\pi N^2\delta}\sum_{k=1}^N\bigg(\sum_{j\ge k+3}^N\frac{1}{j-k-1}+\sum_{j\le k-3}\frac{1}{k-j-1}\bigg)
\le \frac{1}{4\pi N^2\delta}\sum_{k=1}^N\sum_{j=2}^N\frac{1}{j}
\le \frac{\log N}{4\pi N\delta}.
\end{align*}
Collecting all terms, we see that the expression under \eqref{f21} is bounded from above by $\frac{3+\log N}{4\pi N\delta}$.
Hence
\begin{align}\label{f25}
\cp(C_{L,\delta})\ge \frac{4\pi N\delta}{3+\log N}\ge \frac{4\pi L}{3(3+\log N)}\ge L \bigg(\log \frac{L}{2\delta}\bigg)^{-1},
\end{align}
where we have used that $N\ge 3, \delta\le L/6 $.
Numerical evaluation gives that $k_3\ge 0.0101...\ge \frac{1}{100}$. The lower bound in Lemma \ref{lem3} follows by \eqref{f18} and \eqref{f25}.
\hspace*{\fill }$\square $

\subsection{Proof of Theorem \ref{the1}}\label{3.3}
We partition $\Omega_L$ into $N^m$ disjoint open cubes $C_1,\dots, C_{N^m}$ each with measure $(L/N)^m$. We denote the centres of these cubes by $c_1,\dots, c_{N^m}$ respectively. Let $0<\delta<\frac{L}{2N}$, and put
\begin{equation}\label{e26}
\O_{\delta,N,L}=\Omega_L-\cup_{i=1}^{N^m} B(c_i;\delta).
\end{equation}
\begin{figure}[h]
\centering\includegraphics[scale=.6]{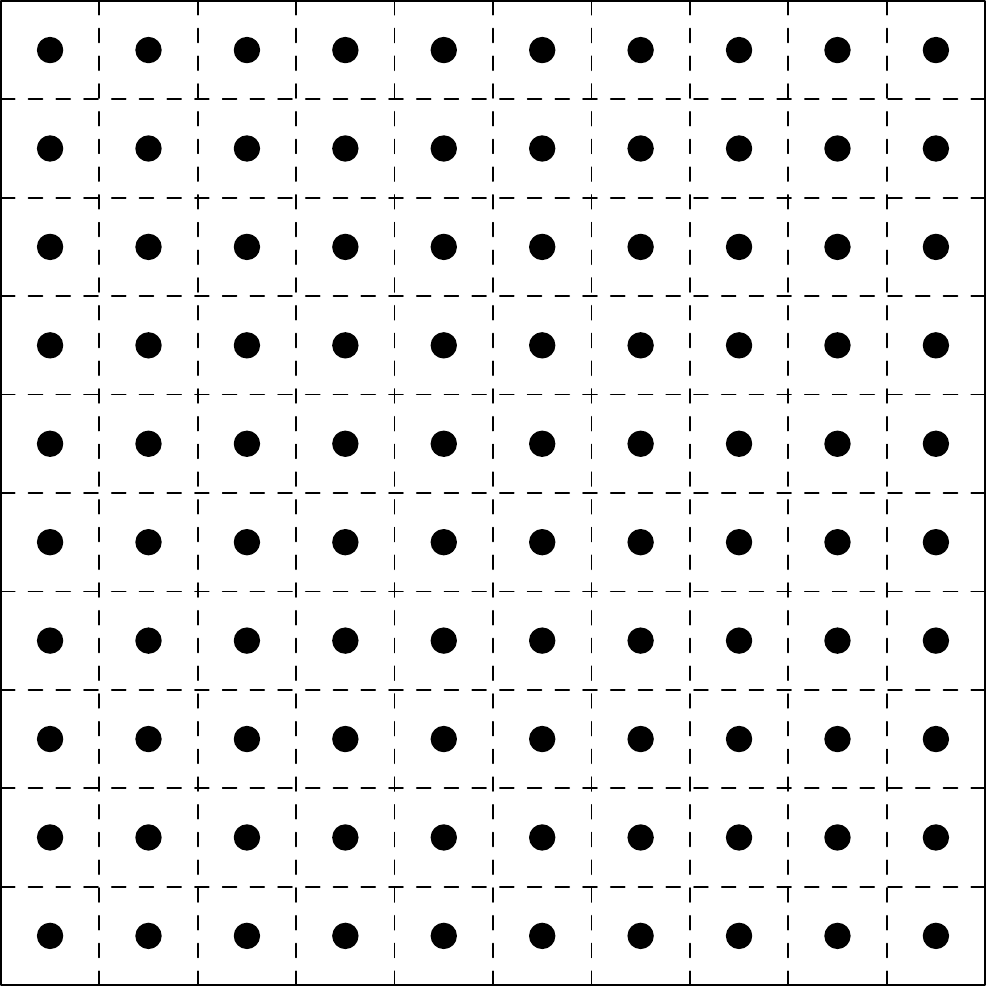}
\caption{$\O_{\delta,N,L}$ with $m=2,N=10,\delta=\frac{L}{8N}.$}\label{fig1}
\end{figure}

 We denote the Dirichlet heat kernel for $\O_{\delta,N,L}$ and $\O_L$ by $p_{\O_{\delta,N,L}}(x,y;t)$ and $p_{\O_L}(x,y;t)$ respectively. The heat kernel with Neumann boundary conditions on $\partial \O_L$ and Dirichlet boundary conditions on $\partial \O_{\delta,N,L}-\partial\O_L$ will be denoted  by
$\pi_{\O_{\delta,N,L}}(x,y;t)$. Let $T>0$, and let $\epsilon>0$ be arbitrary. We bound the torsional rigidity for $\O_{\delta,N,L}$ from below as follows.
\begin{align}\label{e27}
T(\O_{\delta,N,L})&=\int_{\O_{\delta,N,L}}dx\,\int_{\O_{\delta,N,L}}dy\, \int_0^{\infty}dt\,p_{\O_{\delta,N,L}}(x,y;t)\nonumber \\ &\ge
\int_{\O_{\delta,N,L}}dx\,\int_{\O_{\delta,N,L}}dy\, \int_0^{T}dt\,p_{\O_{\delta,N,L}}(x,y;t)\nonumber \\ &=\int_{\O_{\delta,N,L}}dx\,\int_{\O_{\delta,N,L}}dy\, \int_0^{T}dt\,\pi_{\O_{\delta,N,L}}(x,y;t)\nonumber \\ &
\hspace{20mm}-\int_{\O_{\delta,N,L}}dx\,\int_{\O_{\delta,N,L}}dy\, \int_0^{T}dt\,(\pi_{\O_{\delta,N,L}}(x,y;t)-p_{\O_{\delta,N,L}}(x,y;t)).
\end{align}
We now use \eqref{e45f}, \eqref{e45g} and \eqref{e45j} with $K=\cup_{i=1}^{N^m} B(c_i;\delta),$ and $A=\O_{\delta,N,L},\, A=\partial \O_{\delta,N,L}-\partial \O_L$ respectively. So
\begin{equation*}%\label{e30}
\P_x[\tau_{\partial \O_{\delta,N,L}}>t]=\int_{\O_{\delta,N,L}}dy\, p_{\O_{\delta,N,L}}(x,y;t),
\end{equation*}
\begin{equation*}%\label{e31}
\tilde{\P}_x[\tilde{\tau}_{\partial \O_{\delta,N,L}-\partial \O_L}>t]=\int_{\O_{\delta,N,L}}dy\, \pi_{\O_{\delta,N,L}}(x,y;t),
\end{equation*}
and
\begin{align*}%\label{e32}
\tilde{\P}_x[\tilde{\tau}_{\partial \O_{\delta,N,L}-\partial \O_L}>t]&=\tilde{\P}_x[\tilde{\tau}_{\partial \O_{\delta,N,L}}>t]+\tilde{\P}_x[\tilde{\tau}_{\partial \O_L}<t<\tilde{\tau}_{\partial\O_{\delta,N,L}}]\nonumber \\ &\le \P_x[\tau_{\partial \O_{\delta,N,L}}>t]+\tilde{\P}_x[\tilde{\tau}_{\partial \O_L}<t]\nonumber \\ &
=\P_x[\tau_{\partial \O_{\delta,N,L}}>t]+\P_x[\tau_{\partial \O_L}<t].
\end{align*}
Hence the second term in the right-hand side of \eqref{e27} is bounded in absolute value by
\begin{align}\label{e33}
\int_{\O_{\delta,N,L}}dx\int_0^T dt\,\P_x[\tau_{\partial \O_L}<t] &\le\int_{\O_L}dx\int_0^T dt\,\P_x[\tau_{\partial \O_L}<t]\nonumber \\ &\le
\int_{\O_L}dx\int_0^T dt\,\P_x[\tau_{\partial B(x;d(x))}<t]\nonumber \\ &\le 2^{(m+2)/2}\int_0^T dt\int_{\O_L}dx\,e^{-d(x)^2/(8t)}\nonumber \\ &
\le 2^{(m+2)/2}\mathcal{H}^{m-1}(\partial \O_L)\int_0^T dt\int_0^{\infty}dr\,e^{-r^2/(8t)}\nonumber \\ &=s_mL^{m-1}T^{3/2},
\end{align}
with
\begin{equation*}%\label{e34}
s_m=2^{(m+7)/2}m\pi^{1/2}/3.
\end{equation*}
In \eqref{e33} we denoted by $d(x)=\min\{\vert x-y\vert:y\in \partial\O_L\}$, and by $\mathcal{H}^{m-1}(\partial \O_L)$ the surface area of $\partial \O_L$.
In the third inequality of \eqref{e33} we used Corollary \cite[6.4]{MvdBD}, while the fourth inequality follows from the fact that parallel sets of a convex set have decreasing surface area. See  \cite[Proposition 2.4.3]{BB}.

By the periodicity of the cooling obstacles in $\O_{\delta,N,L}$ and the fact that we have no heat flow across $\partial \O_L$ we conclude that
\begin{equation*}%\label{e35}
\int_{\O_{\delta,N,L}}dx\,\int_{\O_{\delta,N,L}}dy\, \,\pi_{\O_{\delta,N,L}}(x,y;t)=N^m\int_{C_1-B(c_1;\delta)}dx\,\int_{C_1-B(c_1;\delta)}dy\,\pi_{B(c_1;\delta), C_1}(x,y;t),
\end{equation*}
where $\pi_{B(c_1;\delta), C_1}(x,y;t)$ denotes the heat kernel with Neumann boundary conditions on $\partial C_1$ and Dirichlet boundary conditions on $\partial B(c_1;\delta)$.
We denote the spectral resolution of the corresponding Laplace operator acting in $\Leb^2(C_1-B(c_1;\delta))$ by $\{\mu_{j,B(c_1;\delta),L/N},j=1,2,\dots\}$, and denote a corresponding orthonormal basis of eigenfunctions by $\{\varphi_{j,B(c_1;\delta),L/N},j=1,2,\dots\}$.
Using the spectral resolution as in \eqref{e19} we have that
\begin{align*}%\label{e36}
\int_{C_1-B(c_1;\delta)}dx\,\int_{C_1-B(c_1;\delta)}dy\,\pi_{B(c_1;\delta), C_1}(x,y;t)&=\sum_{j\in\N} e^{-t\mu_{j,B(c_1;\delta),L/N}} \left(\int_{C_1-B(c_1;\delta)}\varphi_{j,B(c_1;\delta),L/N}\right)^2\nonumber \\ &\ge e^{-t\mu_{1,B(c_1;\delta),L/N}}\|\varphi_{1,B(c_1;\delta),L/N}\|^2_{\Leb^1(C_1-B(c_1;\delta))}.
\end{align*}
We conclude that the first term in the left-hand side of \eqref{e27} is bounded from below by
\begin{align}\label{e37}
N^m&\mu_{1,B(c_1;\delta),L/N}^{-1}(1-e^{-T\mu_{1,B(c_1;\delta),L/N}} )\|\varphi_{1,B(c_1;\delta),L/N}\|^2_{\Leb^1(C_1-B(c_1;\delta))}\nonumber \\ &\ge N^m\mu_{1,B(c_1;\delta),L/N}^{-1}\|\varphi_{1,B(c_1;\delta),L/N}\|^2_{\Leb^1(C_1-B(c_1;\delta))}-\mu_{1,B(c_1;\delta),L/N}^{-1}L^me^{-T\mu_{1,B(c_1;\delta),L/N}},
\end{align}
where we have used \eqref{e38}.
By \eqref{e27}, \eqref{e33} and \eqref{e37} we have that
\begin{align}\label{e39}
T(\O_{\delta,N,L})\ge N^m\mu_{1,B(c_1;\delta),L/N}^{-1}&\|\varphi_{1,B(c_1;\delta),L/N}\|^2_{\Leb^1(C_1-B(c_1;\delta))}\nonumber \\&-\mu_{1,B(c_1;\delta),L/N}^{-1}L^m
e^{-T\mu_{1,B(c_1;\delta),L/N}}-s_mL^{m-1}T^{3/2}.
\end{align}
By Dirichlet-Neumann bracketing (\cite{RS}), we first replace the Dirichlet boundary conditions on $\partial\O_L$ by Neumann boundary conditions and we subsequently insert Neumann boundary conditions on all of the boundaries of the cubes $C_i$. This gives that $\lambda_1(\O_{\delta,N,L})\ge \mu_{1,B(c_1;\delta),L/N}$. Furthermore $\vert \O_{\delta,N,L}\vert\le L^m$. So by \eqref{e4}, \eqref{e39} and Lemma \ref{lem1}, we obtain that
\begin{align}\label{e40}
F(\O_{\delta,N,L})&\ge \frac{N^m}{L^m}\|\varphi_{1,B(c_1;\delta),L/N}\|^2_{\Leb^1(C_1-B(c_1;\delta))}-e^{-T\mu_{1,B(c_1;\delta),L/N}}-s_m\mu_{1,B(c_1;\delta),L/N}L^{-1}T^{3/2}\nonumber \\ &
\ge 1-\left(\frac{4mL^2\mu_{1,B(c_1;\delta),L/N}}{3N^2e}\right)^{1/2}-e^{-T\mu_{1,B(c_1;\delta),L/N}}-s_m\mu_{1,B(c_1;\delta),L/N}L^{-1}T^{3/2}.
\end{align}
We now choose
\begin{equation}\label{e41}
T=\frac{\log N}{\mu_{1,B(c_1;\delta),L/N}}.
\end{equation}
This gives by \eqref{e40}, \eqref{e41} that
\begin{equation}\label{e42}
F(\O_{\delta,N,L})\ge 1-\left(\frac{4mL^2\mu_{1,B(c_1;\delta),L/N}}{3eN^2}\right)^{1/2}-N^{-1}-s_m\frac{(\log N)^{3/2}}{L\mu_{1,B(c_1;\delta),L/N}^{1/2}}.
\end{equation}
We first consider the case $m\ge 3$, and use the bound in \eqref{e45c} to obtain that for
\begin{equation}\label{e420}
\kappa_m\bigg(\frac{N\delta}{L} \bigg)^{m-2}\le \frac{1}{16},
\end{equation}
\begin{equation}\label{e4200}
\bigg(\frac{4mL^2\mu_{1,B(c_1;\delta),L/N}}{3eN^2}\bigg)^{1/2}\le \bigg(\frac{8\pi m^2\kappa_m }{3e}\bigg)^{1/2}\bigg(\frac{N\delta}{L}\bigg)^{(m-2)/2}.
\end{equation}
Similarly \eqref{e45a} gives that
\begin{equation}\label{e4201}
\frac{L}{N}\mu_{1,B(c_1;\delta),L/N}^{1/2}\ge (k_m\kappa_m)^{1/2}\bigg(\frac{N\delta}{L}\bigg)^{(m-2)/2}.
\end{equation}
Combining \eqref{e42}, \eqref{e4200} and \eqref{e4201}, we see that the right-hand side of \eqref{e42} is of the form
\begin{equation}\label{e421}
F(\O_{\delta,N,L})\ge 1-N^{-1}-\bigg(\frac{8\pi m^2\kappa_m }{3e}\bigg)^{1/2}\theta-\frac{s_m}{(k_m\kappa_m)^{1/2}}\frac{(\log N)^{3/2}}{N\theta}.
\end{equation}
with
\begin{equation*}%\label{e422}
\theta=\bigg(\frac{N\delta}{L}  \bigg)^{(m-2)/2}.
\end{equation*}
The choice $\theta=\frac{(\log N)^{3/4}}{N^{1/2}}$ gives that
\begin{equation}\label{e424}
F(\O_{\delta,N,L})\ge 1-O\bigg(\frac{(\log N)^{3/4}}{N^{1/2}}\bigg),\,  N\ra \infty.
\end{equation}

For $m=2$ we obtain by \eqref{e42} and Lemma \ref{lem3} that
\begin{equation}\label{e426}
F(\O_{\delta,N,L})\ge 1-N^{-1}-\bigg(\frac{32}{3e(4-\pi)}\bigg)^{1/2}\theta-\frac{300s_2(\log N)^{3/2}}{\pi N}\theta^{-1},
\end{equation}
with
\begin{equation*}%\label{e427}
\theta=\bigg(\log \frac{L}{2N\delta}\bigg)^{-1/2}.
\end{equation*}
Maximising the right-hand side of \eqref{e426} with respect to $\theta$ yields again \eqref{e424} after a straightforward calculation.
The assertion in Theorem \ref{the1} follows by taking $\Omega_{\epsilon}=C_{\delta,N,L}$ where $\delta, N$ and $L$ satisfy the above relations (for $m=2$ and $m\ge 3$) for the optimal choice of $\theta$, and by choosing $N\in\N$ so large that the term $(\log N)^{3/4}/N^{1/2}$ in \eqref{e424} is smaller than $\epsilon$.
 \hspace*{\fill }$\square $

\section{Proofs of Theorem \ref{the2} and \ref{the3}}\label{sec3}
In this section we give the proofs of Theorem \ref{the2} and Theorem \ref{the3} respectively.

\noindent{\it Proof of Theorem \ref{the2}.}
By the John's ellipsoid Theorem (\cite{J}), there exists an ellipsoid $\Upsilon$ with centre $c$ such that
$\Upsilon \subset \Omega \subset c + m(\Upsilon - c).$
Here $c + m(\Upsilon - c)=\{c + m(x - c) : x \subset \Upsilon\}.$ This is
the dilation of $\Upsilon$ by a factor of $m$ with centre $c$. $\Upsilon$ is the ellipsoid of maximal volume in $\Omega$.
By translating both $\Omega$ and $\Upsilon$ we may assume that
\begin{equation*}
\Upsilon=\{x\in \R^m: \sum_{i=1}^m \frac{x_i^2}{a_i^2}<1\}, \qquad a_i>0,\quad i=1,\dots,m.
\end{equation*}
It is easily verified that the unique solution of \eqref{e1} for $\Upsilon$ is given by
\begin{equation*}
v_{\Upsilon}(x)=2^{-1}\left(\sum_{i=1}^m\frac{1}{a_i^2}\right)^{-1}\left(1-\sum_{i=1}^m\frac{x_i^2}{a_i^2}\right).
\end{equation*}
By changing to spherical coordinates, we find that
\begin{equation}\label{k1}
T(\Omega)\ge T(\Upsilon)=\int_{\Upsilon} v_{\Upsilon}=\frac{\omega_m}{m+2}\left(\sum_{i=1}^m\frac{1}{a_i^2}\right)^{-1}\prod_{i=1}^ma_i.
\end{equation}
Since $\O\subset m\Upsilon$,
\begin{equation}\label{k2}
|\O|\le\int_{m\Upsilon}dx=\omega_mm^m\prod_{i=1}^ma_i.
\end{equation}
By monotonicity of Dirichlet eigenvalues, we have that $\lambda_1(\O)\ge \lambda_1(m\Upsilon)$. The ellipsoid $m\Upsilon$ is contained  in a cuboid with
lengths $2ma_1,\dots,2ma_m.$ So we have that
\begin{equation}\label{k3}
\lambda_1(\O)\ge \frac{\pi^2}{4m^2}\sum_{i=1}^m\frac{1}{a_i^2}.
\end{equation}
Combining \eqref{k1}, \eqref{k2} and \eqref{k3} gives the lower bound in \eqref{e8}.

To prove part (ii) we note (see \cite{EM}) that for bounded, convex $\Omega$ in $\R^2$,
\begin{equation}\label{k4}
\lambda_1(\Omega)\ge \frac{\pi^2\textup{Per}(\Omega)^2}{16\vert\Omega\vert^2}.
\end{equation}
Furthermore, by \cite[Theorem 5.1]{DPG}, we have that for $\Omega$ convex in $\R^m$,
\begin{equation}\label{k5}
T(\Omega)\ge \frac{\vert\Omega\vert^3}{3\textup{Per}(\Omega)^2}.
\end{equation}
The assertion under \eqref{e9} follows by \eqref{k4} and \eqref{k5}.
 \hspace*{\fill }$\square $

\noindent{\it Proof of Theorem \ref{the3}.}

We claim that it is always possible to choose $z_1, z_2\in\partial\Omega$ such that $|z_1-z_2|=w$, and therefore the vector $z_1-z_2$ is orthogonal at $z_1$ and $z_2$ to two parallel supporting hyperplanes achieving the minimal distance $w$.

To show this, the first step is to prove that for any direction $\nu$, there exist two points $\tilde z_1,\tilde z_2\in \partial \Omega$ such that the supporting hyperplanes tangent
to $\partial\Omega$ at these points are parallel to each other. Indeed, assuming that the set is smooth and strictly convex (the general case follows at once from an approximation argument),
for every $\eta\in\mathcal{S}^{m-1}$ such that $\eta\cdot \nu>0$, there exists a unique point $\tilde x(\eta)\in\partial\Omega$ where the outer unit normal is $\eta$. Moreover, there exists
a unique point $\bar x(\tilde x)\in \partial \Omega$ such that $\tilde x-\bar x$ is parallel to $\nu$. We denote by $\xi(\tilde x)$ the inner unit normal to $\Omega$ at $\tilde x$ and observe that $\xi\cdot \nu>0$.
Therefore, denoting by $\mathcal{S}_\nu=\{\eta\in\mathcal S,\eta\cdot \nu\ge0\}$, the map $\xi(\tilde x(\bar x(\eta)))$ (possibly extended so that $\xi=-\eta$ when $\eta\cdot\nu=0$) is a continuous map from $\mathcal{S}_\nu$ into
itself. Brouwer's fixed point theorem provides the existence of $\hat\eta$ such that $\xi(\hat\eta)=\hat\eta$ and this completes the first step.
Now, in view of the above result, %it is always possible to chose $\tilde y$ and $\bar y\in\partial\Omega$ such that $|\tilde y-\bar y|=w$. Indeed
assuming that $T_1$ and $T_2$ are two supporting hyperplanes at distance $w$, there exist two points $z_1,z_2\in \partial \Omega$ such that $z_1-z_2$ is orthogonal to $T_1$ and $T_2$, and the supporting hyperplanes tangent to $\partial\Omega$ at $z_1$ and $z_2$ are parallel to each other. On one hand we have $w\le |z_1-z_2|$, and on the other hand, by construction, $|z_1-z_2|$, is not greater than the distance between $T_1$ and $T_2$. This forces $z_1$ and $z_2$ to belong to $T_1$ and $T_2$ and hence $w= |z_1-z_2|$, which proves our claim.

We introduce a reference frame in $\mathbb{R}^m$, $(x,y)\in\mathbb{R}\times\mathbb{R}^{m-1}$ where $x$ points in the direction $z_1-z_2$ and $(0,0)=\dfrac{z_1+z_2}2$.
Denoting by $E$ the projection of $\Omega$ onto the hyperplane $x=0$, we have
\begin{equation}\label{c0}
\Omega=\{(x,y)\in \mathbb{R}^m: l(y)<x<L(y),\, y\in E\},
\end{equation}
where $L:E\mapsto \mathbb{R}$ is concave,
$l:E\mapsto \mathbb{R}$ is convex, $l\le L$ and $\max\{L(y)-l(y):\, y\in E \}=w$. This maximum is achieved at $y=0$.

We note that $\{(x,y)\in \R^m: x=0, (x,y)\in \O\}\supset \frac12 E$, where $\frac12 E$ is the homothety of $E$ by $\frac12$ with respect to $y=0$. We consider the two-sided cone with base $\frac12 E$ and vertices $(\frac{w}{2},0)$ and  $(\frac{-w}{2},0)$. Let $h\in [0,\frac{w}{2}]$. This two-sided cone contains a cylinder $C_h$ with height $2h$ and base  $\big (1-\frac{2h}{w} \big) \frac12 E$. By monotonicity of Dirichlet eigenvalues, we have that $\lambda_1(\O)\le \lambda_1(C_h)$. By separation of variables, we have that
\begin{equation}\label{c1}
\lambda_1(\O)\le \frac{\pi^2}{4h^2}+\frac{4w^2 \Lambda}{(w-2h)^2}.
\end{equation}
Minimising the right-hand side of \eqref{c1} with respect to $h$ gives that
\begin{align}\label{c2}
\lambda_1(\O)\le \frac{\pi^2}{w^2}\bigg(1+\frac{3c}{2}+\frac{3c^2}{4}+\frac{c^3}{8}\bigg),
\end{align}
where $c$ is given by \eqref{c22}.

If we denote the torsion function of $\O$ by $v=v(x,y)$ where $x\in \R$ and $y\in \R^{m-1}$, then
\begin{align}
\frac{|\Omega| }{T(\Omega)}& =  |\Omega| \frac {\displaystyle\int_E \left(\int_{l}^{L}( |D_y v|^2+v_x^2)\, dx\right) \, dy}  {\left(\displaystyle\int_E\left( \int_{l}^{L} v\,dx\right) dy\right)^2 } \ge  |\Omega| \frac {\displaystyle\int_E \left(\int_{l}^{L} v_x^2\, dx\right) \, dy}  {\left(\displaystyle\int_E\left( \int_{l}^{L} v \,dx\right) dy\right)^2 }\nonumber   \\
& \geq \inf_\phi |\Omega| \frac {\displaystyle\int_E \left(\int_{l}^{L} \phi_x^2\, dx\right) \, dy}  {\left(\displaystyle\int_E\left( \int_{l}^{L} \phi \,dx\right) dy\right)^2 } = \frac {\displaystyle\int_E (L -l)\, dy}{\displaystyle\int_E \frac {(L-l)^3}{12}\, dy}.
\end{align}
The last equality follows by the fact that the function
\begin{equation*}%\label{c5}
(x,y)\mapsto \frac{1}{2}\left\{\left(\frac{L(y) - l(y)} {2}\right)^2 -  \left(x-\frac{L(y) + l(y)} {2}\right)^2\right\}
\end{equation*}
achieves the minimum.
We conclude that
\begin{equation}\label{c6}
\frac{T(\Omega)}{ |\Omega| } \le \frac{1}{12} \frac {\displaystyle\int_E(L -l)^3\, dy}{\displaystyle\int_E(L-l)\, dy}
\le \frac{w^2}{12}.
\end{equation}
Combining \eqref{c2} with \eqref{c6} gives \eqref{a2}.

To prove part (ii), we note that for $m=2$ Theorem \ref{the4} gives that for any $\O$ with finite Lebesgue measure,
\begin{equation*}%\label{c7}
F(\O)\le 1- \frac{\pi }{\lambda_1(\O)|\O|+\pi}.
\end{equation*}
By Blaschke's theorem, \cite[p.215]{Y}, $\O$ contains a ball with radius $w/3.$ Hence $\lambda_1(\O)\le 9{j_{0,1}}^2/w^2$, where $j_{0,1}=2.405\dots$ is the first positive zero of the Bessel function $J_0$. Furthermore, since $|\O|\le w|E|$ and $|E|\ge w$, we have that
\begin{equation}\label{c8}
F(\O)\le 1-\frac{\pi}{\pi+9j_{0,1}^2}\left(\frac{w}{|E|}\right).
\end{equation}
For $\frac{w}{|E|}$ small we use part (i) to obtain an upper bound. Since
\begin{equation}\label{c26}
\Lambda =\frac{\pi^2}{|E|^2}\le\frac{\pi^2}{w^2} ,
\end{equation}
we have, by \eqref{c22}, that $c\le (32)^{1/3}$. By \eqref{a2} and \eqref{c26}, we get that
\begin{align}\label{c27}
F(\O)\le &\frac{\pi^2}{12}\left(1+\left(\frac{3}{2}+\frac{3}{2^{1/3}}+2^{1/3}\right)c \right)\nonumber \\ &
=\frac{\pi^2}{12}+\big(2^{-4/3}+2^{-2/3}+3^{-1}\big)\pi^2\left(\frac{w}{|E|}\right)^{2/3}.
\end{align}

For $\frac{w}{|E|}$ small we use \eqref{c27} as an upper bound, while for $\frac{w}{|E|}$ large we use \eqref{c8} as an upper bound.
The cross-over point value of $\frac{w}{|E|}$ where the right-hand side of \eqref{c8} equals the right-hand side of \eqref{c27} is bounded from below by
$0.0015197$. This, together with the bound under \eqref{c8}, gives the assertion under \eqref{a3}.

 \hspace*{\fill }$\square $

Below we list some known numerical values of $F$ for some convex planar shapes.

\begin{center}
\renewcommand{\arraystretch}{2}
\begin{tabular}{ |c|c| }
 \hline
 Shape & $F(\textup{Shape})$  \\
 \hline
 Rectangle with sides $a,b$ & $\frac{\pi^2}{12}(1+O(a/b)),\,a/b\downarrow 0$, ($\frac{\pi^2}{12}\approx 0.822$)  \\
 \hline
 Disc & $\frac{j_{0,1}^2}{8}\approx 0.723$ \\
 \hline
   Half-disc & $\big(\frac14-\frac{2}{\pi^2}\big)j_{1,1}^2\approx 0.695$  \\
  \hline
Equilateral triangle & $\frac{\pi^2}{15}\approx 0.658$\\
\hline
\end{tabular}
\end{center}

In the table above $j_{1,1}$ is the first positive zero of the Bessel function $J_1$. The values for the thin rectangle and the disc are taken from \cite{MvdBBV}. The torsional rigidity of an equilateral triangle with side lengths $a$ equals $\frac{\sqrt 3a^4}{320}$ \cite[pp. 263--265]{TG}. To obtain the third line in the table we note that for a half-disc of area $\pi a^2/2$, the torsional rigidity is given by $\big(\frac{\pi}{8}-\frac{1}{\pi}\big)a^4$ (see
\cite[pp. 265--267]{TG}). The first Dirichlet eigenvalue of the half-disc is the second Dirichlet eigenvalue of the full disc and equals $j_{1,1}^2a^2$. The first Dirichlet eigenvalue of an equilateral triangle with side lengths $a$ is given by $\frac{16\pi^2}{3a^2}$. So we obtain the last line in the table above.

\end{document}